\documentclass[12pt]{article}
\usepackage{amssymb,amsmath}
\newtheorem{thm}{Theorem}
\newtheorem{cor}[thm]{Corollary}
\newtheorem{lemma}[thm]{Lemma}
\newenvironment{defin}{\medskip\noindent{\sc
Definition}.}{\goodbreak\medskip}
\newenvironment{nota}{\medskip\noindent{\sc
Notations}.}{\goodbreak\medskip}
\newenvironment{remk}{\noindent{\sc
Remark}.}{\goodbreak\vskip10pt}
\newtheorem{prop}[thm]{Proposition}

\def\demo{\medskip\goodbreak\noindent
     \hbox{\sc Proof \kern .3em}\ignorespaces}%
  \def \qedbox{$\square$}%
  \def \qed{\hglue1mm\hfill{\ifmmode\qedbox
     \else\unskip\ \hglue0mm\hfill\qedbox\medskip
      \goodbreak\fi}}%
\def\enddemo{\qed\goodbreak\vskip10pt}%

\def\qed{\hglue1mm\hfill\raise -2pt\hbox{\vrule\vbox to 10pt{\hrule width
4pt
                  \vfill\hrule}\vrule}}

\newcommand{\T}{\mathbb {T}}

\newcommand{\R}{\mathbb {R}}

\newcommand{\Ee}{\mathbb {L}_-}
\newcommand{\Ff}{\mathbb {L}_+}
\newcommand{\Ll}{\mathbb {L}}

\newcommand{\N}{\mathbb {N}}
\newcommand{\Nc}{\mathcal {N}}
\newcommand{\Oc}{\mathcal {O}}
\newcommand{\Vc}{\mathcal {V}}
\newcommand{\Uc}{\mathcal {U}}
\newcommand{\Pc}{\mathcal {P}}

\newcommand{\Cc}{\mathcal {C}}
\newcommand{\Ic}{\mathcal {I}}
\newcommand{\Rc}{\mathcal {R}}
\newcommand{\Mc}{\mathcal {M}}
\newcommand{\Kc}{\mathcal {K}}

\newcommand{\Gc}{\mathcal {G}}
\newcommand{\Lc}{\mathcal {L}}
\newcommand{\Fc}{\mathcal {F}}
\newcommand{\Ac}{\mathcal {A}}

\begin{document}
\title{The tiered Aubry set for autonomous Lagrangian functions}
\author{M.-C. ARNAUD
\thanks{Universit\'e d'Avignon et des Pays de Vaucluse, Laboratoire d'Analyse non lin\' eaire et G\' eom\' etrie (EA 2151),  F-84 018Avignon,
France. e-mail: Marie-Claude.Arnaud@univ-avignon.fr}}
\maketitle
\abstract{ Let $L~: TM \rightarrow {\R}$ be a
 Tonelli Lagrangian function (with $M$ compact and connected and  $\dim M\geq 2$). The tiered Aubry set (resp. Ma\~ n\'e set) $\Ac^T(L)$ (resp. $\Nc^T(L)$) is the union of the Aubry sets (resp. Ma\~ n\'e sets) $\Ac(L+\lambda)$ (resp. $\Nc(L+\lambda)$) for $\lambda$ closed 1-form. Then~:
\begin{enumerate}
\item the set $\Nc^T(L)$ is closed, connected and if $\dim H^1(M)\geq 2$, its intersection with
any energy level is connected and chain transitive;
\item for $L$ generic in the Ma\~ n\'e sense, the sets $\overline{\Ac^T(L)}$ and $\overline{\Nc^T(L)}$ have no interior;
\item if  the interior of $\overline{\Ac^T(L)}$ is non empty, it contains a dense subset of periodic points. 
\end{enumerate}
Then, we give an example of an explicit Tonelli Lagrangian function   satisfying 2  and an example proving that when $M=\T^2$, the closure of the tiered Aubry set and the
closure of the union of the K.A.M. tori may be different.

}

\begin{center}
{\bf R\'esum\'e}
\end{center}
{Soit  $L~: TM \rightarrow {\R}$ un lagrangien de Tonelli (avec $M$ compacte et connexe  et   $\dim M\geq 2$). L'ensemble d'Aubry (resp. de Ma\~ n\'e) \'etag\'e $\Ac^T(L)$ (resp. $\Nc^T(L)$) est la r\'eunion des ensembles d'Aubry (resp. de Ma\~ n\'e) $\Ac(L+\lambda )$ (resp. $\Nc (L+\lambda)$) pour $\lambda$ 1-forme ferm\'ee. On montre~: 
\begin{enumerate}
\item   $\Nc^T(L)$ est ferm\'e, connexe et si $\dim H^1(M)\geq 2$, sa trace avec chaque niveau d'\'energie est connexe et transitive par cha\^\i ne;
\item si $L$ est g\'en\'erique au sens de Ma\~ n\'e, les ensembles  $\overline{\Ac^T(L)}$ et $\overline{\Nc^T(L)}$ sont d'int\'erieur vide;
\item si l'int\'erieur de $\overline{\Ac^T(L)}$ est non vide, il contient une partie dense de points p\'eriodiques.
\end{enumerate}
On donne ensuite un exemple explicite satisfaisant 2 et un exemple montrant que si   $M=\T^2$, $\overline{\Ac^T(L)}$ peut \^etre diff\'erent de l'adh\'erence de la r\'eunion des tores K.A.M. }

\tableofcontents
\newpage
 
\section{Introduction}

 Let $M$ be a compact and  connected  manifold endowed with a
Riemannian metric; we assume that $\dim M\geq 2$. We will denote by
$(x, v)$ a point of the tangent bundle $TM$ with $x\in M$ and $v$ a vector tangent
at
$x$. The projection 
$\pi: TM\rightarrow M$ is then $(x, v)\rightarrow x$. The notation $(x, p)$ will
designate a point of the cotangent bundle $T^*M$ with $p\in T^*_xM$. and $\pi^*: T^*M\rightarrow M$ will be the canonical projection $(x, p)\rightarrow x$.

We consider a Lagrangian  function $L: TM\rightarrow \R$ which is $C^\infty$
and: \begin{enumerate}
\item[$\bullet$] uniformly superlinear: uniformly on $x\in M$, we have: $\displaystyle{\lim_{\| v\|\rightarrow +\infty}\frac{L(x, v)}{\| v\|} =+\infty}$;
\item[$\bullet$] strictly convex: for all $(x, v)\in TM$, $\frac{\partial^2
L}{\partial v^2}(x, v)$ is positive definite.
\end{enumerate} 
Such a Lagrangian function will be called a ``{\em Tonelli Lagrangian
function}''.\\

We can associate to such a Lagrangian  function  the Legendre map
$\Lc=\Lc_L: TM\rightarrow T^*M$ defined by: $\Lc(x, v)=\frac{\partial
L}{\partial v}(x, v)$ which is a fibered
$C^\infty$ diffeomorphism and the   $C^\infty$ Hamiltonian function $H: T^*M\rightarrow
\R$ defined by: $H(x,p)=p\left( \Lc^{-1}(x,p)\right) -L(\Lc^{-1}(x,p))$
(such a Hamiltonian function will be called a ``{\em Tonelli Hamiltonian
function}''). The Hamiltonian function 
$H$ is then superlinear, strictly convex in the fiber and $C^\infty$. We
denote by
$(f^L_t)$ the Euler-Lagrange flow associated to $L$ and $(\Phi_t^H)$ the
Hamiltonian flow associated to $H$; then we have~: $\Phi_t^H=\Lc\circ
f_t^L\circ\Lc^{-1}$.\\

If $\lambda$ is a ($C^\infty$) closed 1-form of $M$, then the map
$T_\lambda~: T^*M\rightarrow T^*M$ defined by~: $T_\lambda (x, p)=(x,
p+\lambda (x))$ is a symplectic ($C^\infty$) diffeomorphism; therefore,
we have~: $(\Phi^{H\circ T_\lambda}_t)=(T_\lambda^{-1}\circ \Phi_t\circ
T_\lambda )$, i.e.  the Hamiltonian flow of $H$ and $H\circ T_\lambda$ are
conjugated. Moreover, the Hamiltonian function $H\circ T_\lambda$ is
associated to the Tonelli Lagrangian function $L-\lambda$, and it is
well-known that~: $(f_t^L)=(f_t^{L-\lambda})$~: the two Euler-Lagrange
flows are equal. Let us emphasize that these flows are equal, but the
Lagrangian functions, and then the Lagrangian actions differ and so  the
minimizing ``objects'' may be different.\\

The reader will find the whole necessary mathematical background
concerning Mather set, Aubry set and Ma\~ n\'e set in the section \ref{S3}.

For   a Tonelli Lagrangian function ($L$ or $L-\lambda$), J.~Mather introduced in \cite{Mat3}
(see  \cite{M2} too) a particular subset $\Ac(L-\lambda)$ of $TM$
which he called the ``static set'' and which is now usually called the
``{\em Aubry set}'' (this name is due to A.~Fathi)\footnote{ These sets extend the notion of
``Aubry-Mather'' sets for the twist maps.}. There exist  different but equivalent definitions of
this set (see
\cite{CIPP} ,  
\cite{F3},
\cite{M2} and section \ref{S3}) and it is known that two closed 1-forms
which are in  the same cohomological class define the same Aubry set~:
$$[\lambda_1]=[\lambda_2]\in H^1(M)\Rightarrow \Ac (L-\lambda_1)=\Ac
(L-\lambda_2).$$
It allows us to   introduce the following notation~: if $w\in H^1(M)$ is
a cohomological class, 
$\Ac_w(L)=\Ac (L-\lambda)$ where $\lambda$ is any closed 1-form belonging
to $w$. $A_w(L)$ is compact, non empty and invariant under $(f_t^L)$.
Moreover, J.~Mather proved in \cite{Mat3} that it is  a Lipschitz graph  above a part of the
zero-section (see
\cite{F3} or section \ref{S3} too).

As we are as interested in the Hamiltonian dynamics as in the
Lagrangian ones, let us  define the dual Aubry set~:

\begin{enumerate}
\item[--] if $H$ is the Hamiltonian function associated to the Tonelli
Lagrangian function $L$, its {\em dual Aubry set} is $\Ac^*(H)=\Lc_L
(\Ac (L))$; 
\item[--] if $w\in H^1(M)$ is a cohomological class, then
$\Ac^*_w(H)=\Lc_L(\Ac_w(L))$ is the {\em $w$-dual Aubry set}; let us
notice that for any closed 1-form
$\lambda$ belonging to $w$, we have~: $T_\lambda(
\Ac^*(H\circ T_\lambda))=\Ac_w^*(H)$.
\end{enumerate}
These sets are invariant under the Hamiltonian flow $(\Phi_t^H)$. 

Another important invariant subset in the theory of Tonelli Lagrangian
functions is the so-called Mather set.  For it, there exists one
definition (which is in     \cite{F3},
\cite{M2}, \cite{Mat1})~: it is the closure of the union of the supports
of the minimizing measures for $L$; it is denoted by $\Mc(L)$ and the {\em
dual Mather set} is $\Mc^*(H)=\Lc_L(\Mc (L))$ which is compact,
non empty and invariant under the flow $(\Phi_t^H)$. As for the Aubry set,
if
$w\in H^1(M)$ is a cohomological class, we define~: $\Mc_w(L)=\Mc
(L-\lambda)$ which is independent of the choice of the closed 1-form
$\lambda$ belonging to $w$. Then
$\Mc^*_w(H)=\Lc_L(\Mc_w(L))=T_\lambda(\Mc^*(H\circ T_\lambda))$ is
invariant under $(\Phi^H_t)$; we name it the {\em $w$-dual Mather set}. 

In a similar way, if  $\Nc (L)$ is  the Ma\~ n\'e set, the {\em dual
Ma\~ n\'e set} is 
$\Nc^*(H)=\Lc_L(\Nc(L))$; we note that if $w\in H^1(M)$ and $\lambda\in
w$, then $\Nc_w(L)=\Nc(L-\lambda)$ is independent  of the choice of
$\lambda\in w$ and then the {\em $w$-dual Ma\~ n\'e set} is
$\Nc_w^*(H)=\Lc_L(\Nc_w(L))=T_\lambda(\Nc^*(H\circ T_\lambda))$; it  is
invariant under $(\Phi_t^H)$, compact and non empty but is not
necessarily a graph.

For every cohomological class $w\in H^1(M)$, we have the inclusion~:
$\Mc^*_w(H)\subset \Ac^*_w(H)\subset  \Nc^*_w(H) $. Moreover, there
exists a real number denoted by
$\alpha_H (w)$ such that~: $\Nc^*_w(H)\subset H^{-1} (\alpha_H (w))$ (see
\cite{Car} and \cite{Mat1}), i.e. each dual Ma\~ n\'e set is contained in an energy level.
For $w=0$, the value $\alpha_H(0)$ is denoted by $c(L)$ and is named
the ``critical value'' of $L$.\\

\begin{defin} If $H~: T^*M\rightarrow \R$ is a Tonelli Hamiltonian
function, the {\em tiered Aubry set}, the {\em tiered Mather set}  and
the {\em tiered Ma\~ n\'e set} are~:
$$\Ac^T(L)=\bigcup_{w\in H^1(M)}\Ac_w(L);\quad \Mc^T(L)=\bigcup_{w\in
H^1(M)}\Mc_w(L);\quad \Nc^T(L)=\bigcup_{w\in
H^1(M)}\Nc_w(L).$$
Their dual sets are~:
$$\Ac^T_*(H)=\bigcup_{w\in H^1(M)}\Ac^*_w(H);\quad
\Mc^T_*(H)=\bigcup_{w\in H^1(M)}\Mc^*_w(H);\quad \Nc^T_*(H)=\bigcup_{w\in
H^1(M)}\Nc^*_w(H).$$
\end{defin}
We shall prove in proposition \ref{P9} that the map  $w\rightarrow \Nc^*_w(L)$ is upper
semi-continuous (roughly speaking, these sets are ``minimizing
objects'' ), therefore   $\Nc_*^T(H)$ is a closed subset of
$T^*M$. It is unknown if such a result
is true or false for the Aubry sets (see
\cite{Mat2}). Concerning the topological structure of the tiered Ma\~ n\'e set, we have~:

\begin{prop}\label{P1} Let $H~: T^*M\rightarrow \R$ be a Tonelli Hamiltonian function.  
Then $\Nc^T_*(H)$    is closed, connected and if $\dim H^1(M)\geq 2$,  for every $h\in \R$, the set
$\Nc^T_*(H)\cap H^{-1}(h)$ is compact, connected and the restriction of $(\Phi_t^H)$ to
$\Nc^T_*(H)\cap H^{-1}(h)$ is chain transitive.
\end{prop}

\noindent{\bf Examples~:} {\bf 1)} At first, let us consider the most
simple completely integrable Hamiltonian function~: $M=\T^n$ and
$H(x,p)=\frac{1}{2}\| p\|^2$. In other words, we consider the geodesic
flow on the flat torus. Then we have~:
$$\forall w\in H^1(M), \Mc^*_w(H)=\Ac^*_w(H)=\Nc^*_w(H)=\{ (x, p);
p=p_0\}$$
where $p_0$ is a constant 1-form; i.e. each of these sets is an invariant
Lagrangian torus, and all these sets fill up the phase space
$T^*\T^n=\Ac^T_*(H)=\Mc^T_*(H)=\Nc^T_*(H)$.

\noindent{\bf 2)} If we perturb a completely integrable Hamiltonian
system for the $C^\infty$ topology, we know that many invariant
tori  will persist (theorems K.A.M.)~: they are dual Mather, Aubry and
Ma\~ n\'e sets for certain cohomological classes. The weak K.A.M. theorems
(see \cite{F3}) give an answer to the following question~: what did 
happen  to the invariant tori which disappeared during the perturbation?
They prove the existence of positively invariant graphs above the zero
section (which are not continuous, but in a certain sense Lagrangian),
each of these graphs containing one dual Aubry set $\Ac_c^*(H)$ (which is
invariant by the Hamiltonian flow $(\Phi_t^H)$) and possibly some pieces
of the stable manifold of this Aubry set.\\

 For the unperturbed system
(completely integrable), we have shown that the Aubry sets fill up the
phase space; but we shall prove that this situation si not generic (the
definition of genericity is just after the theorem)~: 
\begin{thm}\label{T1} Let $H~: T^*M\rightarrow \R$ be a $C^\infty$ generic
Tonelli Hamiltonian function\footnote{Let us recall that we assume that $\dim M\geq 2$; the result
is false if $\dim M=1$.}. Then there exists a dense
$G_\delta$ subset
$G(H)$ of $\R$ such that for every $h\in G(H)$, then
$\overline{\Ac^T_*(H)}\cap H^{-1}(h)$ has no interior in $H^{-1}(h)$;

\noindent in particular, the interior of $\overline{\Ac^T_*(H)}$ is
empty.
\end{thm}
In 1996, R.~Ma\~n\'e introduced the notion of ``generic Lagrangian 
function '' in
\cite{M1}:\\
 {\em ``a certain property holds for a generic Lagrangian $L$ if, given a
strictly convex and superlinear Lagrangian
$L_0$, there exists a residual subset $\Oc\subset C^\infty (M)$ such that
the given property holds for every Lagrangian $L$ of the form $L=L_0+\psi$,
$\psi\in\Oc$''.}\\

Then we define (it is the dual definition for the Hamiltonian
functions)~:

\begin{defin} A certain property holds for a generic Hamiltonian $H$ if,
given a Tonelli Hamiltonian function
$H_0$, there exists a residual subset $\Oc\subset C^\infty (M)$ such that
the given property holds for every Hamiltonian $H$ of the form
$H=H_0-\psi$,
$\psi\in\Oc$.
\end{defin}
The theorem \ref{T1} proves that in the sense of the Baire's category the
tiered Aubry set is small (it is not true for the measure category
when $M=\T^n$~: see the (strong) K.A.M. theorems). We may ask ourselves the same
question for the Ma\~ n\'e set. Let us recall~:

\begin{defin} If $A$ is a closed set invariant under the flow
$(\Phi_t^H)$, its stable manifold $W^s(A, (\Phi^H_t))$ (resp. its
unstable manifold
$W^u(A, (\Phi^H_t))$) is defined by~:
$$W^s(A, (\Phi^H_t))=\{ \xi\in T^*M; \lim_{t\rightarrow +\infty}
d(\Phi^H_t(\xi), A)=0\}$$
(resp~:
$$W^u(A, (\Phi^H_t))=\{ \xi\in T^*M; \lim_{t\rightarrow -\infty}
d(\Phi^H_t(\xi), A)=0\}).$$
\end{defin}
Then it is known  that $\Nc^*_c(H)\subset
W^s(\Ac^*_c(H); (\Phi_t^H))\cap W^u(\Ac^*_c(H); (\Phi_t^H))$ (see \cite{F3} for example); we obtain~:

\begin{cor}\label{C2}
For $H~: T^*M\rightarrow \R$ generic, the set $W^s\left (\overline{\Ac_*^T(H)}; (\Phi_t^H)\right)\cup
W^u\left(\overline{\Ac_*^T(H)}; (\Phi_t^H)\right)$ has no interior; in particular,
the tiered Ma\~ n\'e set
$\Nc^T_*(H)$ has no interior.
\end{cor}

\begin{remk}
A usual default of the genericity results is that in general, one proves these results
by using Baire's theorem but one cannot exhibit one single example of a ``generic object''. It is
not the case of our result, and we obtain easily examples of such ``generic Tonelli Hamiltonian
functions''. Let us consider a Hamiltonian function $H$ whose flow is Anosov on every regular
energy level (for example the geodesic flow on a surface with negative curvature). Then the
restriction of the Hamiltonian flow to every connected regular level is transitive; therefore, if
$c$ is the critical value of $H$,  the set
$D$ of the points
$p$ of $T^*M$ whose orbit is dense in their energy level $E_p$ above $c$  is a dense $G_\delta$
subset
$G$ of
$\{ x\in TM; H(x)>c\}$. But it is known that every orbit of a point of $\Nc^T_*(H)$ is a
lipschitz graph above a part of the zero section, and therefore doesn't meet $G$~: its
interior is then empty.

\end{remk}

Before proving   theorem \ref{T1}, we shall prove the two following results, the first one
explaining in particular what happens in  the ``non generic case'' (when the interior of
$\overline{\Ac^T(L)}$ is non empty), the second one stating precisely which generic property
we need~:

\begin{prop}\label{Per}
Let $L~: TM\rightarrow \R$ be a Tonelli Lagrangian function; then the interior of
$\overline{\Ac^T(L)}$ has a dense subset of periodic points with period in $\N^*$ whose
orbits   are graphs  above a part of the zero section, and have no conjugate point.
\end{prop}

\begin{prop}\label{Pgene} For the Tonelli Lagrangian functions, the
following property is generic~: 
``if $P$ is a periodic orbit which is a graph above a part of the
zero section, which has no conjugate point and  whose period is an
integer $N\geq 1$, then~:
\begin{enumerate}
\item[$\bullet$] either $P$ is hyperbolic and isolated among the
$N$-periodic orbits;
\item[$\bullet$] or in every neighborhood of $P$ there exists an open
subset of points whose orbit has    conjugate points''.
\end{enumerate}
\end{prop}

As we know that every orbit contained in $\overline{\Ac^T(L)}$ has no conjugate point (see
section \ref{S2} for   details), the last  assertion of theorem \ref{T1} is an easy
consequence  of propositions
\ref{Per} and \ref{Pgene}. Then, we use a standard Baire argument to deduce the first
assertion of theorem
\ref{T1} (see the end of section \ref{S4} for the details).

Let us now mention and comment on some related
 results (in these results, the notion of genericity is not exactly the one which we defined
before)~:
\begin{enumerate}
\item in \cite{He2},  M.~Herman proved that for $C^\infty$ generic
exact symplectic twist maps of the annulus $T^*\T^1$, there doesn't exist 
any invariant curve containing a periodic point (section I.5)~: it
implies that for such a generic twist map, the closure of the union of
the Aubry-Mather sets has no interior;  
\item in \cite{He1}, M.~Herman announced~:   for
the
$C^\infty$ generic exact symplectic diffeomorphisms $f$ of $T^*\T^n$
which are homotopic to the identity, if $I(f)$ is the closure  of the union     the  invariant
K.A.M. tori, then $I(f)$  has no interior;
\item I proved in \cite{Arna1} that the generic $C^\infty$ symplectic
diffeomorphisms of any compact symplectic manifold verify~: ``the
closure of the union of the periodic K.A.M. tori in equal to the closure
of the set of the completely elliptic periodic points''; this implies of course
that for such a generic symplectic diffeomorphism, the closure of the
union of the K.A.M. tori has no interior.
\end{enumerate}
We may ask ourselves the following question~: when $M=\T^n$, is our set (the closure of
$\Ac^T(L)$) different from the one introduced by M.~Herman in \cite{He1}? Let us give a
definition~:

\begin{defin} Let $H~: T^*\T^n\rightarrow \R$ be a Tonelli hamiltonian function. A K.A.M. torus
for $H$ is   a Lagrangian $C^\infty$ graph $G$ above
the zero section which is invariant by $(\Phi_t^H)$ and such that the
restriction of the flow to $G$  is conjugated to an ergodic flow~:
$(\theta\rightarrow \theta+t\alpha)$. \\
The closure of the union of the K.A.M. tori for $H$ is denoted by $\Ic(H)$.
\end{defin}
These K.A.M. tori are in fact
Lipschitz graphs and on every compact, the Lipschitz constant may be
chosen uniformly (it is only an adjustment of the results of
\cite{He1}). Moreover, it is proved in \cite{F3} that every exact Lagrangian invariant $C^1$ graph
is the graph of the derivative of a so called ``weak K.A.M. solution'', and therefore every  
Lagrangian invariant $C^1$ graph meets an Aubry set (for a certain cohomology class); if such
a Lagrangian graph is a K.A.M. tori, it is in fact   a Aubry set (because its dynamic
are minimal). We deduce~:
$\Ic (H)\subset
\Nc^T_*(H)$.

We shall build an example such that the tiered Aubry set is not in the closure of the union of the
K.A.M. tori~:

\begin{prop}\label{PCEX} There exists a $C^\infty$ Tonelli Lagrangian function $L~:
T\T^2\rightarrow
\R$ and an open subset $\Uc$ of $C^\infty (\T^2)$ which contains $0$ such that~: for every $\psi\in
\Uc$, there exists a periodic orbit $P$ for $(f_t^{L+\psi})$ which belongs to $\Ac^T(L+\psi)$ ,
but $\Lc_L(P)$ is not in the closure of the union of the  K.A.M. tori for $H-\psi$ ($H$ is the
Hamiltonian function associated to $L$).

\begin{remk}
This result is not very surprising; it corresponds to the existence of Birkhoff instability
regions for twist maps~: in these regions, there exist  periodic orbits which are some Mather
sets, but there exists no ``K.A.M. curve''.
\end{remk}

\end{prop}

\section{  Peierls barrier,   Ma\~ n\'e   potential, Aubry and
Ma\~ n\'e sets, proof of proposition  \ref{P1}}\label{S2}

We gather in this sections some well-known results; the ones concerning
the Peierls barrier are essentially due to A.~Fathi (see
\cite{F3}), the others concerning Ma\~ n\'e potential are given in
\cite{M1}, \cite{CDI} and \cite{CI}.\\
At the end of this section, we prove some new results and proposition 1.

 In the whole section,
$L$ is a Tonelli Lagrangian function.
At first, let us introduce some
notations (we simplify the notation of the critical value~:
$c=c(L)$)~:

\begin{nota}\begin{enumerate}
\item[$\bullet$] given two points $x$ and $y$ in $M$ and $T>0$, we
denote by $\Cc_T(x,y)$ the set of absolutely continuous curves $\gamma~:
[0, T]\rightarrow M$ with $\gamma (0)=x$ and $\gamma (T)=y$;
\item[$\bullet$]  the Lagrangian action along an absolutely continuous
curve
$\gamma~: [a, b]\rightarrow M$ is defined by~:
$$A_L(\gamma)=\int_a^bL(\gamma(t), \dot\gamma (t))dt;$$
\item[$\bullet$] for each $t>0$, we define the function $h_t~:
M\times M\rightarrow \R$ by~: 
$h_t(x,y)=\inf\{ A_{L+c}(\gamma);
\gamma\in
 \Cc_t(x,y)\}$;
\item[$\bullet$] the Peierls barrier is then the function $h~: M\times
M\rightarrow \R$ defined by~:
$$h(x,y)=\liminf_{t\rightarrow + \infty} h_t(x,y);$$ 
\item[$\bullet$]   we define the {\em (Ma\~ n\'e)  
potential} $m~: M\times M\rightarrow \R$ by~: $m 
(x,y)=\inf\{ A_{L+c}(\gamma);
\gamma\in
\bigcup_{T>0}\Cc_T(x,y)\}=\inf\{ h_t(x,y); t>0\}$.
\end{enumerate}
\end{nota}
Then, the Ma\~ n\'e potential verifies~: 
\begin{prop} We have~:
\begin{enumerate}
\item $m$ is finite and $m\leq h$;
\item $\forall x, y, z\in M, m(x,z)\leq m(x,y)+m(y,z)$;
\item $\forall x\in M, m(x,x)=0$;
\item if $x,y\in M$, then $m(x,y)+m(y,x)\geq 0$;
\item if $M_1=\sup\{ L(x,v); \| v\| \leq 1\}$, then~: $\forall x, y\in
M, |m(x,y)|\leq (M_1+c)d(x,y)$;
\item $m~: M\times M\rightarrow \R$ is $(M_1+c)$-Lipschitz. 
\end{enumerate}
\end{prop} 
Now we can define~: 

\begin{defin}\begin{enumerate}
\item[$\bullet$] a absolutely continuous curve $\gamma~:
I\rightarrow M$ defined on an interval $I$ is a {\em ray} if~:
$$\forall [a, b]\subset I,
A_{L+c}(\gamma_{|[a,b]})=h_{(b-a)}(\gamma(a), \gamma(b));$$
a ray is always a solution
of the Euler-Lagrange equations;
\item[$\bullet$] a absolutely continuous curve $\gamma~:
I\rightarrow M$ defined on an interval $I$ is {\em semistatic} if~:
$$\forall [a, b]\subset I, m_c(\gamma (a), \gamma
(b))=A_{L+c}(\gamma_{|[a,b]});$$ a semistatic curve is always a ray;
\item[$\bullet$] the {\em Ma\~ n\'e set} is then~: $\Nc(L)=\{ v\in TM;
\gamma_v\quad is \quad semistatic\}$ where $\gamma_v$ designates the
solution $\gamma_v~: \R\rightarrow M$  of the Euler-Lagrange equations with initial condition  $v$
for
$t=0$; $\Nc (L)$  is  contained in the critical energy level;
\item[$\bullet$] a absolutely continuous curve $\gamma~:
I\rightarrow M$ defined on an interval $I$ is {\em  static} if~:
$$\forall [a, b]\subset I, -m_c(\gamma (b), \gamma
(a))=A_{L+c}(\gamma_{|[a,b]});$$ a  static curve is always a semistatic
curve;
\item[$\bullet$] the {\em Aubry set} is then~: $\Ac(L)=\{ v\in TM;
\gamma_v\quad is \quad  static\}$.

\end{enumerate}
\end{defin}
The following result is proved in \cite{CI}~:
\begin{prop}\label{Pr4}
If $v\in TM$ is such that $\gamma_{v|[a, b]}$ is static for some $a<b$,
then $\gamma_v~: \R\rightarrow M$ is static, i.e. $v\in \Ac(L)$.
\end{prop}
The Peierls barrier verifies (this proposition contains some
results of \cite{CP}, \cite{F3} and \cite{Co1})~:
\begin{prop} (properties of the Peierls barrier $h$) 
\begin{enumerate}
\item the values of the map $h$ are finite  and $m\leq h$;
\item if $M_1=\sup\{ L(x,v); \| v\| \leq 1\}$, then~:
$$\forall x,y,x',y'\in M, |h(x,y)-h(x',y')|\leq
(M_1+c)(d(x,x')+d(y,y'));$$   
therefore $h$ is Lipschitz;
\item  if $x,y\in M$, then $h(x,y)+h(y,x)\geq 0$; we deduce~: $\forall
x\in M, h(x,x)\geq 0$;
\item $\forall x, y, z\in M, h(x,z)\leq h(x,y)+h(y,z)$;
\item \label{I5} $\forall x\in M, \forall y\in \pi(\Ac (L)),
m(x,y)=h(x,y)\quad{\rm and}\quad m(y,x)=h(y,x)$;
\item $\forall x\in M,h(x,x)=0\Longleftrightarrow x\in \pi(\Ac(L))$.
\end{enumerate}
\end{prop}  
 The last item of this
proposition gives us a characterization of the projected Aubry set
$\pi(\Ac(L))$. Moreover,  we have~:
\begin{prop}\label{P4} (A.~Fathi, \cite{F3}, 6.3.3) When $t$ tends to
$+\infty$, uniformly on $M\times M$, the function $h_t$ tends to the
Peierls barrier
$h$.
\end{prop}
A corollary of this result is given in \cite{CI}~:
\begin{cor} (\cite{CI}, 4-10.9) All the rays defined on $\R$ are semistatic.

\end{cor}

Let us give some properties of the Aubry and Ma\~ n\'e sets (see \cite
{M2} and \cite{CDI})~:

\begin{prop}\label{P8} Let $L~: TM\rightarrow \R$ be a Tonelli lagrangian
function. Then~:
\begin{enumerate}
\item[$\bullet$] the Aubry and Ma\~ n\'e set are compact, non empty  and
$\Ac(L)\subset
\Nc (L)$;
\item[$\bullet$] the Aubry set is a Lipschitz graph above a part of the
zero section;
\item[$\bullet$] if $\gamma~: \R\rightarrow M$ is semistatic, then
$(\gamma, \dot\gamma)$ is a Lipschitz graph above a part of the
zero section;
\item[$\bullet$] the $\omega$ and $\alpha$-limit sets of every point of the Ma\~ n\'e set are
contained in the Aubry set.
\end{enumerate}

\end{prop}
We denote by  $\Lambda^1(M)$ the set of ($C^\infty$) closed 1-forms of
$M$ and $\Kc (TM)$ the set of non empty compact subset of $ TM$. Let us
now prove~:

\begin{prop} \label{P9} Let $L~: TM\rightarrow \R$ be a Tonelli
Lagrangian function.  The map $K~: (\psi , \lambda)\in C^\infty (M)\times
\Lambda^1(M)\rightarrow
\Nc (L+\psi+\lambda)\in \Kc (TM)$ is upper semi-continuous.
\end{prop} 

 \noindent\hbox{\sc Proof of   proposition \ref{P9}~:\kern .3em}

Let $H~: T^*M\rightarrow \R$ be the Tonelli Hamiltonian function associated to $L$.  we prove~: 

\begin{lemma}\label{L13}The
map
$(\psi,
\lambda)\in C^\infty (M)\times
\Lambda^1(M)\rightarrow 
\alpha_{H-\psi}(\lambda)=c(L+\psi-\lambda) $ is continuous. 
\end{lemma}

 \noindent\hbox{\sc Proof of   lemma \ref{L13}~:\kern .3em} We use the characterization of the critical
value with the holonomic (probability) measures (see \cite{M2}  or \cite{CI} for the exact
definition of holonomic measure)~:
$-c(L)$ is the minimum of $A_L(\mu)$ among the holonomic measures $\mu$; then  each such minimizing
measures is invariant under $(f_t^L)$ and is contained in the energy level
$\Lc_L^{-1}(H^{-1}(c(L)))$.\\ To prove that $(\psi,
\lambda)\rightarrow 
\alpha_{H-\psi}(\lambda)=c(L+\psi-\lambda) $ is continuous, we only need to prove the continuity
at $(0,0)$.\\
 As $L$ is superlinear, there exists a compact $K\subset TM$ and a neighborhood
$\Vc$ of $(0, 0)$ in $C^\infty(M)\times \Lambda^1(M)$ such that~: for every $(\psi,
\lambda)\in \Vc$, for every holonomic measure $\mu$ such that the  support of $\mu$ meets
$TM\backslash K$,
$\mu$ is not minimizing for $L+\psi-\lambda$. Indeed, let us fix $\mu_0$  any holonomic
measure on $TM$; there exists a neighborhood $\Vc_0$ of $(0, 0)$ in $C^\infty(M)\times
\Lambda^1(M)$ and a constant $\ell\in\R$ such that~: $\forall (\psi,
\lambda)\in \Vc_0, A_{L+\psi-\lambda}(\mu_0)\leq\ell$.  Because $L$ is superlinear, there
exists a constant $C_1\in\R$ such that~: $\forall (\psi,
\lambda)\in \Vc_0, \forall (x, v)\in TM, (L+\psi-\lambda) (x,v)\leq \ell\Rightarrow \|
v\|\leq C_1$. The Hamiltonian function $H$ associated to $L$ being superlinear too,
there exists a constant $C\in \R$ such that~: if $(\psi, \lambda)\in \Vc_0$, if $(x, v)$ and
$(x_0, v_0)$ are in the same energy level for $L+\psi-\lambda$ and if $\| v_0\| \leq C_1$,
then $\| v\|\leq C$.
Hence, if $(\psi,
\lambda)\in \Vc_0$, if
$\mu_{L+\psi-\lambda}$ is a minimizing measure for $L+\psi-\lambda$,  we have~:
$A_{L+\psi-\lambda}(\mu_{L+\psi-\lambda})=-c(L+\psi-\lambda)\leq \ell$. It implies that
there exists $(x,v)\in {\rm supp}(\mu_{L+\psi-\lambda})$ such that $(L+\psi-\lambda)(x,v)\leq
\ell$ and then $\| v\|\leq C_1$. But, $\mu_{L+\psi-\lambda}$ being minimizing, every point of
its support has the same energy as $(x,v)$ and then~: $\forall (X, V)\in {\rm
supp}(\mu_{L+\psi-\lambda}), \| V\|\leq C$. We choose $K=\{ (x, v); \| v\|\leq C\}$.
\\
   We have then to minimize a continuous function $\mu\rightarrow
A_{L+\psi-\lambda}(\mu )$ on a compact set (the set of holonomic probabilities with support in
$K$), we know that the minimum depends continuously on $(\psi, \lambda)$.

\enddemo

From  lemma \ref{L13} and the fact that $\Nc(L)\subset \Lc^{-1}(H^{-1}(c(L)))$, we deduce that
the  Ma\~ n\'e set cannot ``explode''~: for every
$(\psi,
\lambda)\in C^\infty (M)\times
\Lambda^1(M)$, there exists a neighborhood $\Vc$ of $(\psi, \lambda) $ and a compact $K$ of $TM$
such that~: $\forall  (\psi', \lambda')\in \Vc, \Nc(L+\psi'-\lambda')\subset K$.

Let us assume that the proposition \ref{P9} is not true. Then there
exists a sequence $(\psi_n, \lambda_n)$ in $C^\infty (M)\times
\Lambda^1(M)$  which converges to $(\psi , \lambda )$ and a
sequence $(x_n, v_n)\in TM$ converging to $(x , v )$ such
that~:
\begin{enumerate}
\item[$\bullet$] $\forall n, (x_n, v_n)\in \Nc (L+\psi_n+\lambda_n)$;
\item[$\bullet$] $(x , v )\notin
\Nc(L+\psi +\lambda  )$.
\end{enumerate}
As $(x , v  )\notin
\Nc(L+\psi +\lambda )$, the arc $(t\rightarrow
\gamma (t)=\pi\circ 
 f^{L+\psi +\lambda 
}_t(x,v))$ is not a ray for the Lagrangian $L+\psi+\lambda$ and there
exists
$[a, b]\subset
\R$ and
$\varepsilon >0$ such that~:  $$A_{L+\psi +\lambda 
}(\gamma_{ |[a, b]})\geq h^{L+\psi +\lambda 
}_{(b-a)}(\gamma   (a),
\gamma   (b))-(b-a)c(L+\psi+\lambda )+\varepsilon.$$
When $n$ tends to the infinite, if we define $\gamma_{n}(t)=\pi\circ
f_t^{L+\psi_n +\lambda_n }(x_n, v_n)$, then $(\gamma_n, \dot\gamma_n)$
converges uniformly on  any compact interval to $(\gamma, \dot\gamma)$.

We deduce that for $n$ big enough,   we have~:
$$A_{L+\psi_n+\lambda_n}(\gamma_{{n}|[a, b]})\geq
h^{L+\psi_n +\lambda_n 
}_{(b-a)}(\gamma_{ n}(a),
\gamma_{
n}(b))-c(L+\psi_n+\lambda_n)+\frac{\varepsilon}{2};$$
therefore $\gamma_n$ is not a ray, and then is not semistatic. It is a
contradiction.
\enddemo
We deduce a part of proposition \ref{P1}~: 

\begin{cor} The tiered Ma\~ n\'e set  $\displaystyle{\Nc_*^T(H)=\bigcup_{c\in
H^1(M)}\Nc^*_c(H) }$ is closed.

\end{cor}

\begin{cor}\label{C11} Let $L~: TM\rightarrow \R$ be a Tonelli Lagrangian
function and let $(x, v)\in TM$ be such that $(x, v)\notin \Nc^T(L)$.Then there exist~:
\begin{enumerate}
\item[$\bullet$] an open neighborhood $U$ of $(x, v)$ in $TM$;
\item[$\bullet$] an open neighborhood $\Uc$ of $0$ in $C^\infty(M)$;
\end{enumerate}
such that~:
$$\forall \psi\in\Uc, U\cap \Nc^T(L+\psi)=\emptyset.$$
\end{cor}
 \noindent\hbox{\sc Proof of   corollary \ref{C11}~:\kern .3em} We know that $(x,v)\notin \Nc^T(L)$  and
that
$\Nc^T(L)$ is closed; thus there exists a compact neighborhood $K$ of $(x, v)$ in
$TM$ such that $K\cap \Nc^T(L)=\emptyset$.\\
Let $H$ be the Tonelli Hamiltonian function associated to $L$. Then
J.~Mather proved (see \cite{Mat1} or \cite{CI}) that $\alpha_H$ is
convex and  superlinear. Therefore there exists a convex compact subset
$C$ of $H^1(M)$ and a real $R>0$  such that~: 
\begin{enumerate}
\item[$\bullet$] $\forall (x_1, v_1)\in K, H( \Lc_L(x_1,
v_1))<\frac{R}{2}$;
\item[$\bullet$] $\forall w\in \partial C, \alpha_H(w)>R$.
\end{enumerate}
We deduce from the proposition \ref{P9} that if we define~:
$\displaystyle{\Nc_{\partial C}(L)=\bigcup_{w\in \partial C}\Nc_w(L)}$,
then the map
$\psi\in C^\infty (M)\rightarrow \Nc_{\partial C} (L+\psi)$ is upper
semi-continuous. Therefore, the map $\psi\in C^\infty (M)\rightarrow
H\circ\Lc_L (\Nc_{\partial C}(L+\psi))=\alpha_{H-\psi}(\partial C)$ is
upper semi-continuous too. We have~:
\begin{enumerate}
\item[$\bullet$] $\alpha_H(\partial C)\subset ]R, +\infty [$;
\item[$\bullet$] the map $\psi\in C^\infty (M)\rightarrow
 \alpha_{H-\psi}(\partial C)$ is upper semi-continuous.
\end{enumerate}
Then, there exists a neighborhood $\Uc$ of $0$ in $C^\infty (M)$ such
that~: $\forall \psi\in U, \alpha_{H-\psi}(\partial C)\subset ]R,
+\infty[$. Moreover, if $\Uc$ is small enough, we have~: $\forall \psi\in
\Uc, (H-\psi)\circ \Lc_L(K)\subset ]-\infty, \frac{R}{2}[$. These two
facts implies that~:
\begin{enumerate}
\item[$\bullet$] $\forall \psi\in U, \forall w\in H^1(M)\backslash C,
\alpha_{H-\psi}(w)>R$;
\item[$\bullet$] $\forall \psi\in U, \forall (x_1, v_1)\in K,
(H-\psi)(\Lc_L(x, v))<\frac{R}{2}$.
\end{enumerate}
Then, if $\psi\in \Uc$, for every $w\in H^1(M)\backslash C$, for every
$(x_1, v_1)\in K$~: $(x_1, v_1)\notin \Nc_w(L+\psi)$.\\
Moreover, the map $\displaystyle{\psi\in C^\infty (M)\rightarrow \Nc_{ 
C} (L+\psi)=\bigcup_{w\in   C}\Nc_w(L)}$ is upper semi-continuous
(proposition \ref{P9}). There exists a neighborhood
$\Vc$ of $0$ in $C^\infty(M)$ such that~: $\forall \psi\in \Vc, \Nc_{ 
C} (L+\psi)\cap K=\emptyset$. \\
We obtain the conclusion of the corollary with $\Uc\cap \Vc$ and the
interior of $K$.
\enddemo

 \noindent\hbox{\sc End of the proof of   proposition \ref{P1}~:\kern .3em} Let $L$ be the Tonelli
Lagrangian function associated to $H$. We have proved in  proposition \ref{P9} that the   map $ 
  w\in 
H^1(M)\rightarrow
\Nc_w  (L  )\in \Kc (TM)$ is upper semi-continuous. Moreover, we know that each Ma\~ n\'e
set $\Nc_w(L )$ is connected (and chain transitive). We deduce~:

\begin{lemma}\label{L16}
For every arcwise connected subset $C$ of $H^1(M)$, the restriction of $(f_t^L)$
to $\displaystyle{\bigcup_{w\in C}\Nc_w(L)}$ is chain transitive ($\displaystyle{\bigcup_{w\in
C}\Nc_w(L)}$ is therefore  connected when $C$ is closed).
\end{lemma}

\noindent\hbox{\sc Proof of   lemma \ref{L16}~:\kern .3em} Let $x$, $y$ be two points of
$\displaystyle{\bigcup_{w\in C}\Nc_w(L)}$ and $T>0$, $\varepsilon >0$; we want to connect $x$ to
$y$ via a $(\varepsilon, T)$-chain. Let $w_1, w_2\in C$ be such that $x\in \Nc_{w_1}(L)$ and $y\in
\Nc_{w_2}(L)$. The set $C$ being arcwise connected, there exists a continuous arc $w~: [0,
1]\rightarrow C$ such that $w(0)=w_1$ and $w(1)=w_2$. The map $t\in [0, 1] \rightarrow
\Nc_{w(t)}(L)$ is then upper semi-continuous. Therefore for every $t_0\in [0, 1]$, there exists
$\alpha(t_0)>0$ such that~: $\forall t\in [0, 1]\cap ]t_0-\alpha (t_0), t_0+\alpha (t_0)[, \rho (
\Nc_{w(t_0)}(L), \Nc_{w(t)}(L))<\frac{\varepsilon}{3}$ where we define~:
$$\rho (A, B)=\sup\{ d(b, A); b\in B\}.$$ We deduce~:  
$\forall t\in [0, 1]\cap ]t_0-\alpha (t_0), t_0+\alpha (t_0)[, d (
\Nc_{w(t_0)}(L), \Nc_{w(t)}(L))<\frac{\varepsilon}{3}$ where, if $A, B\subset TM$, we
define~: $\displaystyle{d(A,B)=\inf_{a\in A, b\in B} d(a,b)}$.\\
$[0, 1]$ being compact, using a finite covering, we find a finite sequence $t_0=0< \dots < t_N=1$
such that~: $\forall j, d(\Nc_{w(t_j)}, \Nc_{w(t_{j+1})})<\frac{\varepsilon}{3}$. Then we
define a (finite) sequence of points~: 
\begin{enumerate}
\item[$\bullet$] $x_0=x$; $x_{2N+1}=y$;
\item[$\bullet$] for every $j\in \{ 0, \dots , N\}$, $x_{2j}, x_{2j+1}\in \Nc_{w(t_{j})}(L)$;
\item[$\bullet$] for every $j\in \{ 1, \dots , N\}$, $d(x_{2j-1}, x_{2j})<\frac{\varepsilon}{3}$.
\end{enumerate}
Every $x_j$ being in the chain recurrent set of $\displaystyle{\bigcup_{w\in
C}\Nc_w(L)}$ and each $x_{2j}$ being connected to $x_{2j+1}$ by a $(\frac{\varepsilon}{3}, T)$
chain of
$\displaystyle{\bigcup_{w\in C}\Nc_w(L)}$, we
obtain easily a chain passing through
$x, x_1,
\dots, x_{2N+1}=y$.
\enddemo Using   lemma \ref{L16} for $C=H^1(M)$, we deduce that $\Nc^T_*(H)$ is chain transitive
and therefore connected.\\

To deduce the end of the proof of proposition \ref{P1}, we assume that $\dim H^1(M)\geq 2$;
in this case, we notice that if $h\in \R$,
$\alpha_H^{-1}(h)$ is arcwise connected (it is either a convex subset of $H^1(M)$ or the boundary
of a compact convex subset of $H^1(M)$ whose dimension is at least 2, which is homeomorphic to a
connected sphere) and closed. Moreover, we have~: $$\Nc^T_*(H)\cap H^{-1}(h)=\bigcup_{w\in
\alpha_H^{-1}(h)} \Nc^*_w(L).$$
\enddemo

\section{Radially transformed set and Aubry set, proof of proposition \ref{Per}}\label{S3}

\begin{defin}
Let $T>0$; we define~:
\begin{enumerate}
\item[$\bullet$] the set $\Rc_T(L)$ of the $T$-radially
transformed points under $(f_t^L)$ is~:
$$\Rc_T(L)=\{ (x,v)\in TM; \pi(f_T^L(x,v))=x\} ;$$
 its dual set is then $\Rc_T^*(H)=\Lc(\Rc_T(L))$;
\item[$\bullet$] the set $\Pc_T(L)$ is the set of the $T$-periodic
orbits of the Lagrangian flow $(f_t^L)$~:
$$\Pc_T(L)=\{ (x,v)\in TM;  f_T^L(x,v)=(x,v)\} ;$$
its dual set is then $\Pc_T^*(H)=\Lc(\Pc_T(L))$.
\end{enumerate}
\end{defin}
We note that $\Pc_T(L)\subset \Rc_T(L)$ and that if $\lambda$ is a
$C^\infty$ closed 1-form, we have~: $\Pc_T(L-\lambda)=\Pc_T(L)$,
$\Rc_T(L-\lambda)=\Rc_T(L)$. \\
Some of the  radially transformed points which we described before are
minimizing in a certain sense~: 

\begin{prop} \label{Pradferm}(and definition) Let $\lambda$ be a closed
$C^\infty$ 1-form of
$M$. Then for every $x\in M$, the set~:
$$\Gamma_T(L, \lambda; x)=\{ \gamma\in \Cc_T(x,x); \forall \eta\in
\Cc_T(x,x), A_{L-\lambda} (\gamma )\leq A_{L-\lambda}(\eta)\}$$
is non empty and each $\gamma\in \Gamma_T(L,\lambda; x)$ is a solution
of the Euler-Lagrange equations.\\
Moreover, if $\mu$ is a closed 1-form such that $[\mu]=[\lambda]$, then
$\Gamma_T(L,\lambda;x)=\Gamma_T(L,\mu;x)$.\\
This allows us to define for every $w\in H^1(M)$~:
$\Gamma_T(L,w;x)=\Gamma_T(L,\lambda ; x)$ if $[\lambda]=w$ and~:
\begin{enumerate}
\item[$\bullet$] $\Rc_T(L,w;x)=\{ (\gamma(0), \dot\gamma (0)); \gamma\in
\Gamma_T(L,w;x)\}$ ;
\item[$\bullet$] $\displaystyle{\Rc_T(L,w)=\bigcup_{x\in
M}\Rc_T(L,w;x)}$;
\end{enumerate}
the sets $\Rc_T(L,w;x)$ and $\Rc_T(L,w)$ are closed and we have~:
$\Rc_T(L,w) \subset \Rc_T(L)$.
\end{prop}
This proposition is an easy  consequence of Tonelli theorem (see
\cite{F3}).\\
Let us explain how the radially transformed minimizing points allow
us to approximate the Aubry set~: 
\begin{prop}\label{Pradial}  Let $w\in H^1(M)$, $\varepsilon >0$  and let
$L~: TM\rightarrow
\R
$ be a Tonelli Lagrangian function. Then there exists $T_0>0$ such
that~: 
$$\forall T\geq T_0, \forall (x,v)\in \Rc_T(L, w), x\in
\pi(\Ac_w(L))\Longrightarrow d((x,v), \Ac_w(L))\leq \varepsilon.$$

\end{prop}
What this last proposition says is~: the family $(\Rc_T(L, w)\cap
\pi^{-1}(\Ac_w(L)))_{T>0}$ of non-empty compact subsets of $TM$ tends
to $\Ac_w(L)$ (for the Hausdorff topology) when $T$ tends to $+\infty$.
This will be one of the main ingredients of our proof of theorem
\ref{T1}, which will give us some points near the Aubry set but not in the
Aubry set (for generic Lagrangian functions).

\noindent\hbox{\sc Proof of   proposition \ref{Pradial}~:\kern .3em}
Let us assume that the result is not true; then we may find a sequence
$(T_n)_{n\in \N}$ in $\R^*_+$ tending to $+\infty$, a  sequence $(x_n,
v_n)$ of points of $\Rc_{T_n}(L, w)\cap \pi^{-1}(\Ac_w(L))$ such that~: 
 $\forall n\in\N, d((x_n, v_n),
\Ac_w(L))\geq \varepsilon$. \\
Now we use   proposition \ref{P4}~: let $\lambda$ be a closed 1-form
such that $[\lambda]=w$; then we know that if we define $h_t^\lambda~:
M\times M\rightarrow \R$ by $h^\lambda_t(x,y) =\inf\{
A_{L-\lambda+\alpha_H(w)}(\gamma);
\gamma\in
 \Cc_t(x,y)\}$ and $h^\lambda (x,y)=\liminf_{t\rightarrow +\infty}
h^\lambda _t(x,y)$, the functions $h^\lambda_t$ tend uniformly to
$h^\lambda$ when $t$ tends to $+\infty$; moreover, we know that
$h^\lambda$ is Lipschitz and zero at every $(x,x)$ with $x\in \Ac_w(L)$.
If
$\gamma_n$ designates the solution of the Euler-Lagrange equations with
initial condition
$(\gamma_n(0),
\dot\gamma_n(0))=(x_n, v_n)$ we have then~: $h_{T_n}^\lambda
(x_n,x_n)=A_{L-\lambda +\alpha_H(w)}(\gamma_n)$ tends to $0$ when $n$
tends to the infinite.\\
The sequence $(x_n, v_n)$ is bounded (it is a consequence of the
so-called ``a priori compactness lemma'' (see \cite{F3}, corollary
4.3.2)); therefore we may extract a converging subsequence~: we call it
$(x_n, v_n)$ again and $(x_\infty, y_\infty)$ is its limit. Let us
notice that $x_\infty\in \pi(\Ac_w(L))$ because $\Ac_w(L)$ is compact.
Moreover, we have~: $d((x_\infty, v_\infty),
\Ac_w(L))\geq \varepsilon$. \\
Let $\gamma_\infty$ be the solution of the
Euler-lagrange equations such that $(\gamma_\infty(0),
\dot\gamma_\infty(0))=(x_\infty, v_\infty)$. We want to prove that
$\gamma_\infty$ is static~: we shall obtain a contradiction. When $n$ is
big enough,
$\gamma_n(T_n)=\gamma_n(0) $ is close to
$\gamma_\infty(0)$ and $\gamma_n(1)$ is close to $\gamma_\infty (1)$.
Let us fix $\eta>0$; then we define $\Gamma_{n}^{\eta}~: [0,
T_n+2\eta ]\rightarrow M$ by~: 
\begin{enumerate}
\item [$\bullet$] $\Gamma^ \eta_{n |[0, 1]}=\gamma_{\infty|[0,
1]}$;
\item [$\bullet$]  $\Gamma^ \eta_{n |[1, 1+\eta]}$ is a short
geodesic joining $\gamma_\infty (1)$ to $\gamma_n(1)$;
\item[$\bullet$] $\forall t\in [1+\eta, T_n+\eta], \Gamma^\eta_n
(t)=\gamma_n(t-\eta)$;
\item[$\bullet$] $\Gamma^ \eta_{n |[T_n+\eta, T_n+2\eta]}$ is a short
geodesic joining $\gamma_n(T_n)$ to $\gamma_\infty (0)$.
\end{enumerate}
If we choose carefully a sequence $(\eta_n)$ tending to $0$, we have~: 
$$\lim_{n\rightarrow \infty} A_{L-\lambda +
\alpha_H(w)}(\Gamma_n^{\eta_n})=\lim_{n\rightarrow \infty} A_{L-\lambda +
\alpha_H(w)}(\gamma_n)=0.$$
Because the contribution to the action of the two small geodesic arcs
tends to zero (if the $\eta_n$ are well chosen), this implies~:
$$A_{L-\lambda+\alpha_H(w)}(\gamma_{\infty |[0,1]})+m^\lambda(\gamma_\infty (1),
\gamma_\infty (0))\leq 0,$$ where $m^\lambda$ designates Ma\~ n\'e potential
for the Lagrangian function $L-\lambda$. We deduce then from the
definition of Ma\~ n\'e potential that $m^\lambda (\gamma_\infty (0),
\gamma_\infty (1))+ m^\lambda(\gamma_\infty (1),
\gamma_\infty (0))= 0$ and that~:
$A_{L-\lambda+\alpha_H(w)}(\gamma_{|[0,1]})=m^\lambda (\gamma_\infty (0),
\gamma_\infty (1))$. It implies then that
$A_{L-\lambda+\alpha_H(w)}(\gamma_{|[0,1]})=- m^\lambda(\gamma_\infty (1),
\gamma_\infty (0))$. Let us notice that, changing slightly
$\Gamma_n^\eta$, we obtain too~: 
$$\forall [a, b]\subset [0, +\infty[,
A_{L-\lambda+\alpha_H(w)}(\gamma_{|[a,b]})=- m^\lambda(\gamma_\infty
(b),
\gamma_\infty (a));$$
therefore $\gamma_{|[0, +\infty[}$ is static. To conclude, we use the
proposition \ref{Pr4}. \enddemo

To finish this section we give a result   which explains why in
general the radially transformed points are not in a  Ma\~ n\'e set~: in
this case, they would be periodic.
\begin{prop}\label{P7}
Let $L$ be a Tonelli Lagrangian function and let $w\in H^1(M)$ be a
cohomology class; then, for every $T>0$, we have~:
$$\Nc_w(L)\cap \Rc_T(L)\subset \Pc_T(L)\cap \Rc_T(L,
w)\cap \Ac_w(L).$$
\end{prop}
\noindent\hbox{\sc Proof of   proposition \ref{P7}~:\kern .3em}
By proposition \ref{P8}, we know that if $(x, v)\in \Nc_w(L)$, then its
orbit is a Lipschitz graph above a part of the zero section.  Therefore,
if $(x, v)\in \Nc_w(L)\cap \Rc_T(L)$, then the orbit of
$(x,v)$ under $(f_t^L)$ is a graph (above a part of the zero section);
as $\pi (f_T(x,v))=\pi (x, v)$, we deduce that~: $f_T(x,v)=(x,v)$~: $(x,
v)$ is $T$-periodic for $(f_t)$, i.e. $(x,v)\in \Pc_T(L)$. Moreover,
$\gamma_v$ is a ray~: therefore it is minimizing between $\gamma_v(0)$
and $\gamma_v(T)$~: $(x, v)\in \Rc_T(L, w)$. \\
We deduce from from proposition \ref{P8} that every periodic orbit contained in $\Nc_w(L)$ is
in $\Ac_w(L)$. Hence $(x, v)\in \Ac_w(L)$.
\enddemo
\begin{cor} Let $L$ be  a Tonelli lagrangian function. Let $(x, v)\in
TM$ be such that~:
$(x,v)\in \Rc_T(L)\backslash \Pc_T(L)$ for some $T>0$. Then there exist~:
\begin{enumerate}
\item[$\bullet$] an open neighborhood $U$ of $(x, v)$ in $TM$;
\item[$\bullet$] an open neighborhood $\Uc$ of $0$ in $C^\infty(M)$;
\end{enumerate}
such that~:
$$\forall \psi\in\Uc, U\cap \Nc^T(L+\psi)=\emptyset.$$
\end{cor}
This result is an easy consequence of   proposition \ref{P7} and  
corollary \ref{C11}.

Let us now prove~:

\begin{prop}\label{PB}Let  $U$ be a non empty open subset of of $TM$ and 
 $L~: TM\rightarrow \R$ be a Tonelli Lagrangian function. Then~:
\begin{enumerate}
\item[--] either there exists a non empty subset $U'\subset U$ such that $U'\cap
\Ac^T(L)=\emptyset$;
\item[--] or there exists $N\in\N^*$ and a sequence $(x_n, v_n)$ of different $N$-periodic
points contained in $\Ac^T(L)\cap U$ such that $\displaystyle{\lim_{n\rightarrow \infty}(x_n,
v_n)=(x_0,v_0)}$; moreover the orbit of every $(x_n, v_n)$ is a graph  above a part or the
zero section, and has no conjugate point. 
\end{enumerate}

\end{prop}
Of course, we deduce proposition 4 from this lemma~: if $U\subset
\overline{\Ac^T(L)}$ is a non empty open subset, we have found $(x, v)\in U$ which is
periodic, whose orbit   is a graph  above a part or the
zero section, and has no conjugate point. 

\noindent\hbox{\sc Proof of proposition \ref{PB}~:\kern .3em}

Let us   consider  a  Tonelli Lagrangian function $L$   and let $U$ be a non empty subset of
$TM$. There are two cases~:
\begin{enumerate}
\item either $U\cap \Ac^T(L)=\emptyset$~: we have the first conclusion;
\item or there exists $(x, v)\in U\cap \Ac_w(L)$ for a certain $w\in
H^1(M)$. 
\end{enumerate}
Let us choose $\alpha >0$ such that $\bar B((x,v), \alpha )\subset U$. We know
that $\Ac_w (L)$ is a compact graph above a part of the zero section
such that~: $T_xM\cap \Ac_w(L)=\{ (x, v)\}$. Therefore,  there exists 
$\varepsilon >0$ such that $\forall (x, v_1)\in T_xM, d((x, v_1),
\Ac_w(L))<\alpha \Rightarrow d((x, v_1), (x,v))<\varepsilon$.   
 By   proposition \ref{Pradial},  there exists
$N_0>0$ such that~: 
$$\forall N\geq N_0, \forall (x,v_1)\in \Rc_N(L, w,; x),  d((x,v_1),
\Ac_w(L))\leq \alpha .$$
We deduce that if $N\geq N_0$~: $\Rc_N(L, w; x)\subset U$.\\
Let us recall that the set $\displaystyle{\Rc_N(L, w)=\bigcup_{y\in M}
\Rc_N(L, w;y)}$ is a closed subset of $TM$ (proposition \ref{Pradferm}).
Thus there exists a neighborhood $U_0$ of $x$ in $M$ such that~:
$\forall y\in U_0, \Rc_{N_0}(L, w; y)\subset U$. Another time, we have
two cases~:
\begin{enumerate}
\item either there exist $x_1\in U_0$ and $(x_1, v_1)\in \Rc_{N_0}(L, w;
x_1)$ such that $(x_1, v_1)\notin \Nc^T(L)$; $\Nc^T(L)$ being closed, we
have the conclusion for the set $U'=U\backslash \Nc^T(L)$.
\item or $\displaystyle{\bigcup_{y\in U_0}\Rc_{N_0} (L,w; y)\subset
\Nc^T(L)}$; then by proposition \ref{P7}, $\displaystyle{\bigcup_{y\in
U_0}\Rc_{N_0} (L,w; y)}$ is a union of periodic orbits with period
$N_0$ contained in $\Ac^T(L)$. These orbits are graphs above a part of the zero section, and
have no conjugate point. \enddemo
\end{enumerate}

\section{Green bundles, conjugate points and proofs of of theorem \ref{T1}  and  corollary
\ref{C2}}\label{S4} All the results contained in this section  except the last proposition are
not new.\\ Let us recall some definitions~:

\begin{defin} Let $L$ be a Tonelli Lagrangian function  defined on $TM$
and $(x,v)\in TM$~:
\begin{enumerate}
\item[$\bullet$]   the ``vertical'' at
$(x,v)\in TM$ is the linear subspace
$V(x,v)=\ker D\pi (x,v)$ of $T_{(x,v)}(T M)$; the vertical at
$(x,p)\in T^*M$ is the (Lagrangian) linear subspace
$V(x,p)=\ker D\pi^* (x,p)=D\Lc_L(V(x,v))$ of $T_{(x,p)}(T^*M)$;
\item[$\bullet$]    the orbit of
$(x,v)$ has a conjugate point if there exists $t\not=t'$ such that
$Df^L_{t-t'}(V(f^L_{t'}(x,v)))\cap V(f^L_{t}(x,v))\not=\{ 0\}$; then we say that $t$ and $t'$ are
conjugate (along the orbit); the definition is the same for
$(x,p)$. 
\end{enumerate}

\end{defin}
We recall some results of \cite{CI2}~:

\begin{prop} 
Let $(x,v)=(\gamma (t_0), \dot\gamma (t_0))$ be a point of a ray  $\gamma~: \R\rightarrow M$ for 
$L$; then its orbit has no conjugate point.
\end{prop}

\begin{prop}   Let  $(x,v)$ be a point of $TM$ which is not a fixed
point of the flow $(f_t^L)$ and which has no conjugate point; then there
exists two $(f_t^L)$ invariant $n$-dimensional subbundles of $T(TM)$,
$G_-$ and
$G_+$, named the  Green bundles defined by~: 
$$G_-(x,v)  =\lim_{t\rightarrow
+\infty}Df^L_{-t}(V(\phi_t(x,v)))\quad{\rm and}\quad
 G_+ (x,v)=\lim_{t\rightarrow +\infty}Df^L_{t}(V(f^L_{-t}(x,v))).$$
Moreover, they are transverse to the vertical and  if we define~:
$\Ee(x,v)=D\Lc (x,v) (G_-(x,v))$ and
$\Ff (x,v)=D\Lc (x,v)(G_+(x,v))$, then~: $\Ee (x,v)$ and $\Ff(x,v)$ are
Lagrangian, their sum is contained in the tangent bundle of the energy level of
$\Lc(x,v)$ and their intersection contains the Hamiltonian vector field.

\end{prop}

\begin{prop}\label{Phyp} Let $(x,v)$ be a $T$-periodic point with no conjugate point  of
$(f_t^L)$ which is not a fixed point of the flow. Then, if the dimension of
$G_-(x,v)+G_+(x,v)$ is 2n-1, this orbit is hyperbolic and for every
vector $W\in T_{(x,v)}(TM)\backslash \left( G_-(x,v)\cap G_+(x,v)\right)$
(where
$G_-(x,v)$ and $ G_+(x,v)$ designate  the Green bundles)~:
the family
$(Df_t^L(x,v)(W))_{t\in\R}$ is unbounded.

\end{prop}

Now we shall give a detailed description of the images of the vertical; to do that we need to
introduce some new notions~: 

\begin{defin}
Let $\Ll_1$, $\Ll_2$ be two Lagrangian subspaces of $T_{(x,p)}(T^*M)$ which are
transverse to the vertical. Then the restrictions of $D\pi^*(x,p)~: T_{(x,p)}(T^*M)\rightarrow
T_xM$ to $\Ll_1$ and $\Ll_2$ are two isomorphisms, named $F_1$ and $F_2$. \\
The {\em relative height} between $\Ll_1$ and $\Ll_2$ is then the quadratic form $Q(\Ll_1, \Ll_2)$
defined on $T_xM$ by~:
$$\forall \delta x\in T_xM, Q(\Ll_1, \Ll_2)(\delta x)=\omega (F_1^{-1}(\delta x), F_2^{-1}(\delta
x)).$$
We
say that
$\Ll_2$ is {\em above} $\Ll_1$ if $Q(\Ll_1, \Ll_2)$ is positive (i.e. if its index is
$0$),  that $\Ll_2$ is {\em strictly above} $\Ll_1$ if $\Ll_2$ is above $\Ll_1$
and the dimension of
$\Ll_1\cap
\Ll_2$ is 0   i.e. if $Q(\Ll_1, \Ll_2)$ is positive definite and that $\Ll_2$ is {\em
semi-strictly above} $\Ll_1$ if $\Ll_2$ is above $\Ll_1$ and the dimension of
$\Ll_1\cap
\Ll_2$ is 1   i.e. if $Q(\Ll_1, \Ll_2)$ is positive with nullity  1.
\end{defin}

\begin{remk}
The definition of the height (slightly different because given in a chart) was given in
\cite{Arna2}. 
\end{remk}
Let us recall some results of \cite{CI2} and \cite{Arna2}~:

\begin{prop}\label{Pheight} We define~: $V_t(x,p)=D\Phi_t^H(V(\Phi_{-t}^H(x,p))$. Then~:
\begin{enumerate}
\item[$\bullet$] Let $\Ll$ be a Lagrangian subspace of $T_{(x,p)}(T^*M)$ which is transverse to
$V(x,p)$;  for
$t>0$ small enough~:
$V_t(x,p)$ is strictly above $\Ll$ which is strictly above $V_{-t}(x,p)$ (``small enough'' is
locally uniform in $(x,p)$);
\item[$\bullet$] if $0<t_1<t_2$ and the orbit has no  conjugate point  between $0$
and $t_2$~:
$V_{t_1}(x,p)$ is strictly above $V_{t_2}(x,p)$;
\item[$\bullet$] if $0<t_1<t_2$ and the orbit has no  conjugate point  between $-t_2$ and $0$~:
$V_{-t_2}(x,p)$ is strictly above $V_{-t_1}(x,p)$;
\item[$\bullet$] if $t$, $t'$ are strictly positive and the orbit has no conjugate point between
$-t$ and $t'$, then $V_{t'}(x, p)$ is strictly above $V_{-t}(x,p)$.
\end{enumerate}

\end{prop}
A first consequence is the well-known~:

\begin{cor}\label{Cab}
Let $(x,p)$ be a point having no conjugate point for $H$; then  $\Ff(x,p)$ is above
$\Ee(x,p)$.
\end{cor}

Another  consequence is~: 

\begin{cor}\label{Conj}
Let $H~: T^*M\rightarrow \R$ be a Tonelli Hamiltonian function. Then the subset $\Uc$ of $C^\infty
(M)\times T^*M$ defined by~: $\Uc=\{ (\psi, (x,p))\in C^\infty
(M)\times T^*M$; there are two conjugate points for $H-\psi$  along the orbit of $(x,p)\}$ is open.
\end{cor}
\noindent\hbox{\sc Proof of   corollary \ref{Conj}~:\kern .3em}
Let us assume that there exist some conjugate points along the orbit of $(x,
p)$ for $H$~:  there exist two real numbers $t_1<t_2$ such that
$D\Phi_{t_2-t_1}^H(V(\phi^H_{t_1}(x,p)))\cap V(\Phi^H_{t_2}(x,p))$ contains at least one non zero
vector, named $Y$. To simplify the notations, we may assume that $t_1=0$ and $t_2=T>0$. We may
assume too that $T$ is the smallest $t>0$ such that $0$ and $t$ are conjugate along the orbit of
$(x,p)$. We have~:
\begin{enumerate}
\item[$\bullet$] if $X$ is a non zero vector belonging to $V(x,p)\cap V_{-T}(x,p)$ and if
$Y=D\Phi_T^H(x,p)X$, for
$u>0$, $Z=D\Phi^H_{-u-T}(Y)=D\Phi_{-u}^H(X)\in V_{-T-u}(\Phi^H_{-u}(x,p))\cap
V_{-u}(\Phi^H_{-u}(x,p))$. 
\item[$\bullet$]   for $u>0$ small enough~: $Z\in
V_{-u}(\Phi^H_{-u}(x,p))$, the Lagrangian subspace $V_{-2u}(\Phi^H_{-u}(x,p))$ is strictly
above  $V_{-u}(\Phi^H_{-u}(x,p))$ (it is the third point of proposition \ref{Pheight}); we choose
then
$u$ such that
$u<T$; then we have~:
$V_{-T-u}(\Phi^H_{-u}(x,p))$ is not above $V_{-2u}(\Phi^H_{-u}(x,p))$~: indeed,
$Z$ belongs   to $V_{-T-u}(\Phi^H_{-u}(x,p))$ and to a Lagrangian
subspace, $V_{-u}(\Phi^H_{-u}(x,p))$,   which is ``strictly under''  $V_{-2u}(\Phi^H_{-u}(x,p))$.
\end{enumerate}
Finally, we have found $(x_1, p_1)$ on the orbit of $(x,p)$ and $0<t_1<t_2$ such that
$V_{-t_2}( x_1,p_1) $ is not above $V_{-t_1}(x_1,p_1 )$; this condition is clearly open
and implies the existence of conjugate points (see proposition \ref{Pheight}).\enddemo

Let us now prove proposition \ref{Pgene}~:

\noindent{\bf Proposition \ref{Pgene} }{\em For the Tonelli Lagrangian functions, the
following property is generic~: 
``if $P$ is a periodic orbit which is a graph above a part of the
zero section, which has no conjugate point and  whose period is an
integer $N\geq 1$, then~:
\begin{enumerate}
\item[$\bullet$] either $P$ is hyperbolic and isolated among the
$N$-periodic orbits;
\item[$\bullet$] or in every neighborhood of $P$ there exists an open
subset of points whose orbit has    conjugate points''.
\end{enumerate}
}

\noindent In the proof of proposition \ref{Pgene}, we will prove the following result, which
is the main ingredient of the proof~:

\begin{prop}  Let $L~: TM\rightarrow \R$ be a Tonelli Lagrangian functions and $P$ be a
non hyperbolic periodic orbit of $(f_t^L)$ which is a graph above a part of the zero section,
which has no conjugate point and  whose period is an integer $N\geq 1$, then~:\\
in every neighbourhood of $0$ in $C^\infty(M, \R)$, there exists a function $\psi$ such that
 $P$ is a  periodic orbit for $(f_t^{L+\psi})$ with conjugate points.
\end{prop}

\noindent\hbox{\sc Proof of   proposition  \ref{Pgene}~:\kern .3em} Let $H~: T^*M\rightarrow \R$ be a
Tonelli Hamiltonian function. Let $(U_n)$ be a countable basis of open and
relatively compact subsets of $T^*M$. The subset $\Uc_n$ of $C^\infty (M)$ is the set of $\psi$
such that $U_n$ contains a point whose orbit under $(\Phi^{H-\psi}_t)$  has a conjugate point. We
deduce from corollary \ref{Conj} that $\Uc_n$ is open. Therefore $\Vc_n= \Uc_n\cup
(C^\infty(M)\backslash
\bar \Uc_n)$ is an open and dense subset of the Baire space $C^\infty (M)$ and
$\Gc=\bigcap_{n\in\N} \Vc_n$ is a dense $G_\delta$ subset of $C^\infty (M)$. \\
Let us consider
$\psi\in \Gc$ and let $(x,p)$ be a $N$-periodic point for $(\Phi^{H-\psi}_t)$ whose orbit is a
graph above a part of the zero section, which has no conjugate point. Let us assume that there
exists a neighborhood  $U_n$  of $(x,p)$  which  contains no point whose orbit under
$(\Phi^{H-\psi}_t)$  has a conjugate point. As $\psi\in \Gc$ and $\psi\notin\Uc_n$, we have~:
$\psi\in C^\infty(M)\backslash \bar \Uc_n$.\\
Let us now consider the orbit of $(x,p)$~: as it has no conjugate points, we can define the Green
bundles $\Ee$ and $\Ff$ along this orbit. There are two cases~:\\
1) If these Green bundles are transverse in the
energy level, we use proposition \ref{Phyp}~: $\Ee (x,p)\cap \Ff(x,p)=\R X_H(x,p)$ where $X_H$ is
the Hamiltonian vector field, the orbit is hyperbolic and the eigenvectors of $D\Phi^H_N(x,p)$
associated to the eigenvalue 1 are the vectors of $\R X_H(x,p)$ (because the orbits of the  other
vectors are unbounded); it implies that this orbit is  isolated
among the
$N$ periodic orbits. \\
2) If the Green bundles are not transverse in the energy level, we shall show that we may add
to
$H-\psi=\tilde H$ a small function $\psi_1\in C^\infty (M)$ to create conjugate points along the
orbit of $(x,p)$; it will imply that $\psi\in \overline{ \Uc_n}$, it is a contradiction with 
$\psi\in C^\infty(M)\backslash \bar \Uc_n$.\\

Let us now build such a $\psi_1$. We assume that $(x,p)$ is not a fixed point of the flow (this
case is  simpler that the case which we treat); then there exists $t_0\geq 0$ such that, if we
define~: $\gamma (t)=\pi^*\circ \Phi_t^{\tilde H}(x,p)$, then $\dot\gamma (t_0)\not=0$; we
define~: $x_0=\gamma(t_0)$. We choose
$C^\infty$-coordinates $(x^1, \dots, x^n)$ in a neighborhood $U\subset
M$ of
$x_0$ such that if $U\cap \gamma=\{ \gamma(t);t\in ]t_0-\varepsilon,
t_0+\varepsilon[\}$, then: $\forall t\in ]t_0-\varepsilon,
t_0+\varepsilon[, (x^1, \dots , x^n)(\gamma(t))=(t, 0\dots , 0)$. We work then
in the dual (symplectic) coordinates $(x^1, \dots , x^n, p^1, \dots , p^n)$ on
$T^*U$~: it means that the point with coordinates $(x^1, \dots, x^n, p^1, \dots p^n)$ is 
$\displaystyle{\sum_{k=1}^n p^kdx^k}$.
 We define a function $\psi_1: M\rightarrow \R$ which is: \begin{enumerate}
\item[--] zero on $M\backslash U$;
\item [--] defined in the chart $U$ by: $\displaystyle{\psi_1 (x)= \eta
\left(\sum_{i=1}^n
\left( x^i\right)^2\right) \sum _{j=2}^n \left( x^j\right)^2}$ where $\eta: \R\rightarrow [0, 1]$ is a
$C^\infty$ function which is zero outside $]-\left(
\frac{\varepsilon}{2}\right)^2,
\left( \frac{\varepsilon}{2}\right)^2[$ and strictly positive  in
$]-\left(
\frac{\varepsilon}{2}\right)^2,
\left( \frac{\varepsilon}{2}\right)^2[$.
\end{enumerate}
Then $(x, p)$ has the same (periodic) orbit $\Gamma$ for $(\Phi_t^{\tilde H})$ as for
$(\Phi_t^{\tilde H+\psi_1})$.\\
Let us now assume that the orbit of $(x,p)$ has
no conjugate point  for $(\Phi_t^{\tilde
H+\psi_1})$ (we shall show a contradiction). Then we may define along this orbit the Green bundles
$\Ee^1$ and
$\Ff^1$ (for
$\tilde H+\psi_1$). We shall use~:\\
\begin{lemma}\label{L27} We consider $(x,p)\in \Gamma$,  $\tau>0$ and  $\Ll $ a Lagrangian
subspaces of
$T_{(x,p)}(T^*M)$ transverse to
$V(x,p)$ such that: \begin{enumerate}
\item[a)] for every $t\in [0, \tau]$, $D\phi^{\tilde H }_t(\Ll )$ is transverse to
$V(\phi^{\tilde H}(x,p))$;
\item[b)] for every $t\in [0, \tau]$, $D\phi^{\tilde H+\psi_1}_t(\Ll )$ is transverse to
$V(\phi^{\tilde H}(x,p))$.
\end{enumerate}
Then for every $t\in [0, \tau]$, $D\phi^{\tilde H }_t(\Ll )$ is
above   (semi-strictly above if $t\geq N$) $D\phi^{\tilde H+\psi_1}_t(\Ll )$.

\end{lemma}
\noindent\hbox{\sc Proof of   lemma  \ref{L27}~:\kern .3em}
We begin by proving a version of this lemma for small $t$.\\
We   say that $(\delta x, \delta p):
\R\rightarrow T(T^*M)$ is an infinitesimal solution along the orbit of $(x,p)$ for $(\phi_t)$ if
$(\delta x(t), \delta p(t))\in T_{\phi_t(x,p)}(T^*M)$ and $(\delta
x(t), \delta p(t))=D\phi_t(\delta x(0), \delta p(0))$. Let $(\delta x_1, \delta
p_1)$ (resp. $(\delta x_0, \delta
p_0)$) be an infinitesimal solution for $\tilde H+\psi_1$ (resp. $\tilde H$)
along $\Gamma$. They satisfy the so-called linearized Hamilton equations
(given in coordinates): $$\delta \dot x_1 =\frac{\partial^2 \tilde
H}{\partial x\partial p}\delta
x_1+\frac{\partial^2
\tilde H}{\partial p^2}\delta p_1;\quad \delta
\dot p_1 =-\frac{\partial^2 \tilde H}{\partial x^2}\delta
x_1-\frac{\partial^2 \tilde H}{\partial p\partial x}\delta p_1-\frac{\partial^2 \psi_1}{\partial
x^2}(x)\delta x_1;$$
$$\delta \dot x_0 =\frac{\partial^2 \tilde
H}{\partial x\partial p}\delta
x_0+\frac{\partial^2
\tilde H}{\partial p^2}\delta p_0;\quad 
\delta \dot
p_0 =-\frac{\partial^2 \tilde H}{\partial x^2} \delta x_0-\frac{\partial^2
\tilde H}{\partial p\partial x}\delta p_0.$$
We are interested in some infinitesimal solutions having the same initial
values: $(\delta x_0(0), \delta p_0(0))=(\delta x_1(0), \delta p_1(0))$. We  deduce from
the linearized Hamilton equations that, uniformly for $(x,p)\in \Gamma$ close to $(x_0, p_0)$, if
the two infinitesimal solutions have the same initial values, for $t$ close to $0$: $$(*)\quad 
\delta x_1 (t)=\delta x_0(t)+O(t^2);\quad \delta p_1(t)=
\delta
p_0(t)-t\frac{\partial^2 \psi_1}{\partial x^2}(x)\delta x_1(t) +O(t^2).$$
Let us   assume that we work in a dual chart and that $\Ll$ is a Lagrangian subspace of
$T_{(x,p)}(T^*M)$ transverse to $V(x,p)$~: then $\Ll$ is the graph of a symmetric matrix $S$~:
$\delta p=S\delta x$. In this chart, the coordinates of $\Phi_t^{\tilde H}(x,p)$ are $(x(t),
p(t))$.  For
$t$ small enough
$D\Phi_t^{\tilde H}\Ll=\Ll_t$ and
$D\Phi_t^{\tilde H+\psi_1}\Ll=\Ll_t^1$ are Lagrangian subspace of $T_{\Phi_t^{\tilde
H}(x,p)}(T^*M)$ which are transverse to the vertical~: they are graphs of $S_t$, $S_t^1$. \\
We distinguish two cases (even if they are not exhaustive)~:
\begin{enumerate}
\item[a)] there exists $\alpha>0$ such that the support of $\psi$ doesn't meet $\{  x(t);
t\in ]0; \alpha[\}$; in this case, for every  $t\in [0, \alpha]$, $D\phi^{\tilde H }_t
(\Ll)=D\phi^{\tilde H+\psi_1}_t(\Ll)$;
\item[b)] there exists $\alpha>0$ such that $\{  x(t);
t\in [0; \alpha]\}$ is in the interior of the support of $\psi_1$; let $(\delta x_0,
\delta p_0)$ and $(\delta x_1, \delta p_1)$ be some infinitesimal solutions as before such that~:
$(\delta x_0(0),
\delta p_0(0))=(\delta x_1(0), \delta p_1(0))\in \Ll$; 
we have~: 
$$\delta x_1(t)=\delta x_0(t)+O(t^2);\quad \delta p_0(t)=S_t(x,p)\delta x_0(t);\quad \delta
p_1(t)=S^1_t(x,p)\delta x_1(t)$$
$$\delta p_1(t)=\delta
p_0(t)-t\frac{\partial^2 \psi_1}{\partial x^2}(x)\delta x_1(t) +O(t^2)=(S_t(x,p)-t\frac{\partial^2
\psi_1}{\partial x^2}(x))\delta x_1(t)+O(t^2)$$
We deduce that~: $S^1_t(x,p)=S_t(x,p)-t\frac{\partial^2
\psi_1}{\partial x^2}(x)+O(t^2)$; therefore~: $S_t(x,p)-S^1_t(x,p)$ a symmetric matrix which is
positive with nullity 1  for $t>0$ small enough. It is the matrix (in the chart) of the
relative eight
$Q(D\phi^{\tilde H+\psi_1}_t(\Ll), D\phi^{\tilde H }_t
(\Ll))$; thus, $D\phi^{\tilde H }_t
(\Ll)$ is semi-strictly above  $D\phi^{\tilde H+\psi_1}_t(\Ll)$ for $t$ small enough.
\end{enumerate}
Let us notice that using a limit, we deduce from the case b) that if $\{  x(t);
t\in ]0; \alpha]\}$ is in the interior of the support of $\psi_1$, we have the same conclusion~:
then we have dealt with all the possible cases for $(x,p)$.\\
Now, to prove the lemma for large $t$, we notice that any symplectic  flow preserves the order
between Lagrangian subspaces~: if $\Ll_1$, $\Ll_2$ are such that $\Ll_2$ is above $\Ll_1$ and such
that for every $t\in [0, \tau]$, $D\phi_t(\Ll_1)$ and $D\phi_t(\Ll_2)$ are transverse to the
vertical, then the relative height $Q(D\phi_t(\Ll_1),D\phi_t(\Ll_2))$ has a kernel varying
continuously with
$t$ and whose dimension is constant (it is $\dim (D\phi_T(\Ll_1\cap \Ll_2))=\dim (\Ll_1\cap
\Ll_2)$); therefore its index is constant.

Let us now prove the first part of lemma \ref{L27}. We have proved that there exists $\alpha >0$
such that~: for every $t\in [0, \alpha]$, $D\phi^{\tilde H }_t(\Ll_1)$ is
above     $D\phi^{\tilde H+\psi_1}_t(\Ll_1)$. Let us define~: $\tau_1=\sup \{ R\in [0, \tau]; 
\forall t\in
[0,
R], Q(D\phi^{\tilde H+\psi_1}_t(\Ll_1), D\phi^{\tilde H }_t(\Ll_1))\geq 0\}$. Let us assume
that $\tau_1<\tau$; at first, we notice that, by continuity, the supremum is indeed a maximum~:
$Q(D\phi^{\tilde H+\psi_1}_{\tau_1}(\Ll_1), D\phi^{\tilde H }_{\tau_1}(\Ll_1))\geq 0$, i.e.
$D\phi^{\tilde H }_{\tau_1}(\Ll_1)$ is above $D\phi^{\tilde H+\psi_1}_{\tau_1}(\Ll_1)$; because
the flow preserves the order between Lagrangian subspaces, we deduce that for any $u\in [0,
\tau-\tau_1]$, $D\phi^{\tilde H }_{\tau_1+u}(\Ll_1)$ is above $D\phi^{\tilde H
}_{u}(D\phi^{\tilde H+\psi_1}_{\tau_1}(\Ll_1))$; but for $u>0$ small enough, we have~: 
$D\phi^{\tilde H
}_{u}(D\phi^{\tilde H+\psi_1}_{\tau_1}(\Ll_1))$ is above $D\phi^{\tilde H+\psi_1
}_{\tau_1+u}(\Ll_1)$; therefore, for $u>0$ small enough~: $\tau_1+u$ contradicts the definition of
$\tau_1$.

To obtain the ``semi-strict'' of the lemma, we notice that along a subarc of the orbit of
$(x,p)$ (in the support of $\psi_1$), we find locally a strict inequality between the
Lagrangian subspaces (it is the case b) before).\enddemo
Using this lemma (the vertical is not transverse to itself, but we may use an image of this
vertical), we obtain~: \\
for every $t>N$, $D\Phi_t^{\tilde H  }(V(\Phi^{\tilde H  }_{-t}(x,p)))$ is semi-strictly
above 
$D\Phi_t^{\tilde H +\psi_1}(V(\Phi_{-t}^{\tilde H +\psi_1}(x,p)))$; when $t$ tends to $+\infty$,
we obtain~: $\Ff (x,p)$ is above $\Ff^1(x,p)$. As $\Ff^1(x,p)$ is invariant under $(D\Phi_t^{\tilde
H+\psi_1})$, we deduce that for every $t>0$~: $D\Phi_t^{\tilde H +\psi_1}(\Ff(\Phi_{-t}^{\tilde H
+\psi_1}(x,p)))$ is semi-strictly above $\Ff^1(x,p)$. But by lemma \ref{L27} and the fact that
$\Ff(x,p)$ is invariant under $(D\Phi_t^{\tilde H})$, we know that, for $t>N$,  $\Ff(x,p)$ is
semi-strictly above
$D\Phi_t^{\tilde H +\psi_1}(\Ff(\Phi_{-t}^{\tilde H +\psi_1}(x,p)))$. Finally, $\Ff (x,p)$ is
semi-strictly above $\Ff^1(x,p)$.\\
In a similar way, we obtain that $\Ee^1 (x,p)$ is
semi-strictly above $\Ee (x,p)$. Finally, we have~:
\begin{enumerate}
\item[$\bullet$]$\Ff (x,p)$ is
semi-strictly above $\Ff^1(x,p)$, i.e $Q(\Ff^1(x,p), \Ff (x,p))$ has a 1-dimension kernel and
is positive ;
\item[$\bullet$] $\Ee^1 (x,p)$ is
semi-strictly above $\Ee (x,p)$ i.e. $Q(\Ee (x,p), \Ee^1 (x,p))$ has a 1-dimension kernel and
is positive ;
\item[$\bullet$] $\dim (\Ee(x,p)\cap \Ff(x,p)\geq 2)$ i.e. $Q(\Ee(x,p), \Ff(x,p)$ is positive and
the dimension of its kernel is at least 2.
\end{enumerate}
Therefore~: $$Q(\Ee^1 (x,p), \Ff^1(x,p))=-Q(\Ee (x,p), \Ee^1 (x,p)) +Q(\Ee(x,p),
\Ff(x,p))-Q(\Ff^1(x,p), \Ff (x,p))    $$
is strictly negative at any vector of $\ker Q(\Ee(x,p),
\Ff(x,p))\backslash \ker Q(\Ff^1(x,p), \Ff (x,p))$ and we obtain a contradiction with
corollary \ref{Cab}.\enddemo

 \noindent\hbox{\sc Proof of theorem \ref{T1}~:\kern .3em}\\
We explained in the introduction how we deduce the last assertion of theorem \ref{T1} fom
propositions propositions
\ref{Per} and \ref{Pgene}.\\
 The first part of the theorem is a consequence of the end of the
theorem and a  result of Baire's theory~: \\
we consider a generic Tonelli Lagrangian function.  The union $\Rc^*$  of the regular level of
$H$ is
  a dense open subset of $\R$. We denote  the set of the regular values of $H$ by $V$.\\
Let us consider $h_0\in V$; then there exists a diffeomorphism, $\Phi~:
H^{-1}(h_0)\times ]-\varepsilon, \varepsilon [\rightarrow U\subset T^*M$
such that~: $\forall \eta\in ]-\varepsilon, \varepsilon[,
\Phi(H^{-1}(h_0)\times
\{\eta\} )=H^{-1}(h_0+\eta )$. Then $A=\Phi^{-1}(\overline{\Ac^T_*(H)}\cap U)$ is a closed
subset of $H^{-1}(h_0)\times ]-\varepsilon, \varepsilon[$ which has no
interior. Let $(U_n)$ be a basis of non empty subsets of
$H^{-1}(h_0)$. We define~: $\Fc_n=\{ \eta\in ]-\varepsilon,
\varepsilon[ ;U_n\times
\{
\eta
\}
\subset ( H^{-1}(h_0)\times
\{
\eta\} ) \cap A\}$. As $A$ is closed, $\Fc_n$ is a closed subset of
$]-\varepsilon, \varepsilon[$. Moreover, as $A$ has no interior, $\Fc_n$
has no  interior; therefore $\displaystyle{F=\bigcup_{n\in\N} \Fc_n}$ has
no interior (Baire's theorem) and $G=h_0+(]-\varepsilon,
\varepsilon[\backslash F)$ is a dense $G_\delta$ subset  of
$]h_0-\varepsilon, h_0+\varepsilon[$ such that~: for every $h\in G$,
$H^{-1}(h)\cap \Ac^T_* (H)$ has no interior in $H^{-1}(h)$.\enddemo

\noindent\hbox{\sc Proof of corollary \ref{C2}~:\kern .3em}
We want to prove that the set 
$$W^s(\overline{\Ac_*^T(H)}; (\Phi_t^H))\cup
W^u(\overline{\Ac_*^T(H)}; (\Phi_t^H))$$ has no interior. Let us assume that it is not true
and  let
$U\subset W^s(\overline{\Ac_*^T(H)}; (\Phi_t^H))\cup
W^u(\overline{\Ac_*^T(H)}; (\Phi_t^H))$ be an open and  non empty    subset. Using
theorem
\ref{T1}, we know that the open set $U'=U\backslash \overline{\Ac^T_*(H)}$ is non empty. By  
Poincar\'e recurrence theorem, almost every point in $U'$ (for the volume form associated to
the symplectic form) is positively and negatively recurrent. But a point of $U'$ is in
$$\left( W^s(\overline{\Ac_*^T(H)}; (\Phi_t^H))\cup W^u(\overline{\Ac_*^T(H)};
(\Phi_t^H))\right)
\backslash \overline{\Ac_*^T(H)}; (\Phi_t^H));$$ therefore, either it is not negatively
recurrent or it is not positively recurrent, which contradicts the fact that $U'$ is non
empty.
\enddemo

\section{ Proof of proposition \ref{PCEX}}

We begin by defining a completely integrable Tonelli Hamiltonian function of $T^*\T^2$, whose flow
is the ``product'' of the flow of a pendulum and the geodesic flow of the circle~: if we identify
$T^*\T^2$ with the set $\T\times\R\times\T\times\R$, if the (global) coordinates are $(\theta_1,
p_1, \theta_2, p_2)\in \T\times\R\times\T\times\R$, the Hamiltonian function $H_0$ is defined by~:
$$H_0(\theta_1,
p_1, \theta_2, p_2)=\frac{1}{2}(p_1^2+p_2^2)+\cos (2\pi\theta_1)-\frac{3}{2};$$
then the Hamiltonian flow of $H_0$ is defined by~: $\Phi_t^{H_0}((\theta_1,
p_1, \theta_2, p_2)=(\varphi_t(\theta_1, p_1), \psi_t(\theta_2, p_2))$ where $(\varphi_t)$ is
the flow of the pendulum and $(\psi_t)$ the geodesic flow of $\T$.\\
Let $w$ be the cohomological class of the 1-form $d\theta_2$. Then~: 
$$\Mc_w^*(H_0)=\Ac_w^*(H_0)=\Nc_w^*(H_0)=\{ (0, 0, t, 1); t\in \T\} .$$
If we perturb slightly $H_0$, we may obtain a Hamiltonian function $H_1$ such that~:
\begin{enumerate}
\item $\Nc_w(H_1)=\Nc_w(H_0)$ is a periodic hyperbolic orbit $P$ (in fact, the Ma\~ n\'e set
$\Nc_w(H )$ depends continuously on $H$);
\item the intersections between the stable manifold $W^s(P, (\Phi_t^{H_1}))$ and the unstable
manifold $W^u(P, (\Phi_t^{H_1}))$ are transverse in the energy level $\Sigma= H_1^{-1}(0)$ of $P$;
\item the surface $S =\{ (\theta_1,   p_1, 0, p_2);   
H_1(\theta_1,   p_1, 0, p_2)=0\} $ is transverse to the flow in the $0$ energy level and near the
point $(0,0,0,1)$;
\item    in any
neighborhood $V$ of $(0,0,0,1)$ in $S$, there exists another neighborhood $U$ of $(0,0,0,1)$ in
$S$ such that $U\subset V$ and such that~$\delta U=\gamma_1\cup\gamma_2\cup\gamma_3\cup\gamma_4$
where the $\gamma_i$ are some arcs such that~: $\gamma_1\cup\gamma_3\subset W^s(P,
(\Phi_t^{H_1}))$ and
$\gamma_2\cup\gamma_4\subset W^u(P, (\Phi_t^{H_1}))$ ; to obtain such a result, we only have to
ask that there is a transverse homoclinic intersection on any local branch of
$  W^s(P, (\Phi_t^{H_1}))\cap S$ and $  W^u(P, (\Phi_t^{H_1}))\cap S$~: then we obtain a kind of
canvas by arcs of $  W^s(P, (\Phi_t^{H_1}))\cap S$ and $  W^u(P, (\Phi_t^{H_1}))\cap S$ around
$(0,0,0,1)$ in $S$.  
\end{enumerate}
The situation which we just described is in fact open in the following sense~:\\
{\em there exist  $\varepsilon>0$ and an open subset $\Uc$ of $C^\infty (M)$ containing $0$ 
such that~: for every $\psi\in\Uc$, for every $h\in [-\varepsilon,  \varepsilon]$,
$(H_1+\psi)^{-1}(h)$ contains one periodic orbit $P'=P(\psi, h)$, the orbit of $\xi$,  close to
$P$  and in any neighborhood $V$ of $\xi$ in $S'=S(\psi, h) =\{ (\theta_1,   p_1, 0, p_2);   
(H_1+\psi)(\theta_1,   p_1, 0, p_2)=h\} $, there exists another neighborhood $U$ of $\xi$ in
$S'$ such that $U\subset V$ and such that~$\delta U=\gamma_1\cup\gamma_2\cup\gamma_3\cup\gamma_4$
where the $\gamma_i$ are some arcs such that~: $\gamma_1\cup\gamma_3\subset W^s(P',
(\Phi_t^{H_1+\psi }))$ and
$\gamma_2\cup\gamma_4\subset W^u(P', (\Phi_t^{H_1+\psi}))$.  }

Moreover, there exists a neighborhood $\Vc\subset \Uc$ of $0$ in $C^\infty (M)$ such that, for
every $\psi\in \Vc$, there exists $h=h(\psi )\in ]-\varepsilon, \varepsilon[$ such that
$\Nc_w(H_1+\psi)=P(\psi, h)$ (we have seen that the Ma\~ n\'e set depends continuously on $\psi$, and
an invariant set contained in a neighborhood of a hyperbolic orbit and in an energy level is
necessarily a periodic orbit). Let us prove that for $\psi\in \Vc$, $\Nc_w(H_1+\psi)$ is not a
subset of $\Ic (H_1+\psi)$.\\
Let us consider $\xi\in  P(\psi, h)\cap S(\psi, h)$ and let us assume that there exists  a
sequence of K.A.M. tori $(T_i)_{i\in\N}$ such that   $\displaystyle{\lim_{i\rightarrow
\infty}d(\xi, T_i)=0}$. Being a Lagrangian invariant torus, each
$T_i$ is in an energy level $(H_1+\psi)^{-1}(h_i)$ with, for $i$ large enough~: $h_i\in
]-\varepsilon,
\varepsilon[$ close to $P(\psi, h_i)$. Moreover, $T_i$ is a graph above the zero section.
Therefore, $\{ (\theta_1, p_1, 0, p_2); (\theta_1, p_1, 0, p_2)\in T_i \}$ is  a curve 
$\Gamma_i$ which is a graph above a circle. Moreover, this curve passes very close to the point
$\xi_i$ which is the point of the periodic orbit $P(\psi, h_i)$ which belongs to
$S(\psi, h_i)$.
This curve is then a curve which is traced on
$S(\psi, h_i)$, which has points very close to $\xi_i$ and other points far from $\xi_i$. 
Therefore it cuts the boundary of any sufficiently small neighborhood of $\xi_i$ in $S(\psi,
h_i)$, and then contains some points of $ W^s(P(\psi, h_i), (\Phi_t^{H_1+\psi}))\cup W^s(P(\psi,
h_i), (\Phi_t^{H_1+\psi}))$. This contradicts the fact that the restriction of the flow to any
K.A.M. torus is minimal. \enddemo

\end{document}